\documentclass[12pt,a4paper]{article}

\usepackage[T1]{fontenc}
\usepackage[utf8]{inputenc}
\usepackage[english]{babel}
\usepackage{lmodern,microtype}
\usepackage[left=28mm,right=28mm,top=38mm,bottom=40mm,headheight=16pt,headsep=9mm]{geometry}
\usepackage{setspace}
\usepackage{amsmath,amssymb,amsthm,mathtools,mathrsfs}
\numberwithin{equation}{section}
\newtheoremstyle{plainspaced}{12pt}{12pt}{\itshape}{}{\bfseries}{.}{.6em}{}
\theoremstyle{plainspaced}
\usepackage{aliascnt}
\newtheorem{theorem}{Theorem}[section]
\newaliascnt{proposition}{theorem}
\newtheorem{proposition}[proposition]{Proposition}
\aliascntresetthe{proposition}
\newaliascnt{lemma}{theorem}
\newtheorem{lemma}[lemma]{Lemma}
\aliascntresetthe{lemma}
\newaliascnt{corollary}{theorem}
\newtheorem{corollary}[corollary]{Corollary}
\aliascntresetthe{corollary}
\newtheoremstyle{uprightspaced}{12pt}{12pt}{\normalfont}{}{\bfseries}{.}{.6em}{}
\theoremstyle{uprightspaced}

\newaliascnt{example}{theorem}

\aliascntresetthe{example}

\theoremstyle{remark}
\newaliascnt{remark}{theorem}
\newtheorem{remark}[remark]{Remark}
\aliascntresetthe{remark}
\newtheorem*{remark*}{Remark}

\usepackage{titlesec}
\titleformat{\section}{\normalfont\large\scshape\centering}{\thesection.}{.65em}{}
\titlespacing*{\section}{0pt}{2.2\baselineskip}{.9\baselineskip}
\titleformat{\subsection}{\normalfont\normalsize\bfseries}{\thesubsection.}{.6em}{}
\titlespacing*{\subsection}{0pt}{1.4\baselineskip}{.5\baselineskip}

\usepackage{fancyhdr}
\newcommand{\shorttitle}{Integral Riemann--Roch via Adams operations}
\usepackage{enumitem}
\usepackage{needspace}
\setlist[enumerate]{label=\textup{(\roman*)},leftmargin=*,topsep=.5em,itemsep=.3em}
\setlist[itemize]{leftmargin=*,topsep=.5em,itemsep=.3em}
\usepackage[hidelinks,hyperfootnotes=false]{hyperref}
\usepackage{titling}
\pretitle{\begin{center}\Large\bfseries}
\posttitle{\par\end{center}\vskip .7em}
\preauthor{\begin{center}\large\begin{tabular}[t]{c}}
\postauthor{\end{tabular}\par\end{center}}
\predate{\begin{center}\normalsize}
\postdate{\par\end{center}}

\usepackage{tikz-cd}
\usetikzlibrary{arrows.meta}
\allowdisplaybreaks[1]
\hypersetup{
  pdftitle={Integral Riemann--Roch via Adams operations},
  pdfsubject={Integral Grothendieck--Riemann--Roch in algebraic K-theory},
  pdfkeywords={Integral Riemann--Roch, Grothendieck--Riemann--Roch, Chow groups,
    algebraic K-theory, Adams operations, topological filtration, Chern classes},
  pdfdisplaydoctitle=true,
  bookmarksnumbered=true
}

\usepackage[nameinlink,noabbrev]{cleveref}
\crefname{theorem}{Theorem}{Theorems}
\crefname{proposition}{Proposition}{Propositions}
\crefname{lemma}{Lemma}{Lemmas}
\crefname{remark}{Remark}{Remarks}
\crefname{corollary}{Corollary}{Corollaries}
\crefname{example}{Example}{Examples}
\crefname{section}{Section}{Sections}
\crefname{equation}{}{}
\creflabelformat{equation}{(#2#1#3)}

\newcommand{\CH}{\operatorname{CH}}
\newcommand{\Td}{\operatorname{Td}}
\newcommand{\ch}{\operatorname{ch}}
\newcommand{\cl}{\operatorname{cl}}

\newcommand{\rk}{\operatorname{rk}}
\newcommand{\Z}{\mathbf Z}
\newcommand{\Q}{\mathbf Q}
\newcommand{\OO}{\mathcal O}
\newcommand{\PP}{\mathbf P}
\newcommand{\lcm}{\operatorname{lcm}}
\newcommand{\CT}{\mathfrak{CT}}
\newcommand{\Pcal}{\mathcal P}

\newcommand{\paperauthors}{A. Navarro Garmendia and J. Navarro Garmendia}
\newcommand{\paperemails}{%
\textsc{Alberto Navarro Garmendia}\\
Departamento de Matem\'atica Aplicada, Universidad Polit\'ecnica de Madrid,\\
Madrid, Spain\\
\href{mailto:alberto.navarro@upm.es}{\texttt{alberto.navarro@upm.es}}\\[1em]
\textsc{Jos\'e Navarro Garmendia}\\
Departamento de Matem\'aticas, Universidad de Extremadura,\\
Badajoz, Spain\\
\href{mailto:navarrogarmendia@unex.es}{\texttt{navarrogarmendia@unex.es}}}
\title{Integral Riemann--Roch via Adams operations}
\author{\paperauthors}
\date{}
\hypersetup{pdfauthor={\paperauthors}}

\begin{document}
\begin{center}
{\Large\bfseries \MakeUppercase{\thetitle}\par}
\vspace{1em}{\small\MakeUppercase{\paperauthors}\par}
\end{center}
\thispagestyle{plain}

\begin{abstract}
Let $f:X\to Y$ be a projective morphism between smooth quasi-projective
schemes over an arbitrary field. We consider the defect obtained by clearing
the denominators in the degree-$n$ part of Grothendieck--Riemann--Roch.
For $m=n+\dim X-\dim Y\geq0$, we construct an explicit integer $V_m$,
depending only on $m$, which annihilates this defect; more precisely,
$V_{m-i}$ annihilates it on the $i$-th step of the topological filtration
of $K_0(X)$. 

The proof combines the leading terms of this filtration with
the Adams--Riemann--Roch theorem.
\end{abstract}

\begingroup
\renewcommand{\thefootnote}{}
\footnotetext{\textit{2020 Mathematics Subject Classification.}
Primary 14C40; Secondary 14C17, 19E08, 19L20.\par
\textit{Key words and phrases.}
Integral Grothendieck--Riemann--Roch, Chow groups, algebraic $K$-theory,
Adams operations, topological filtration, Chern classes.}
\endgroup
\setcounter{footnote}{0}

\setcounter{tocdepth}{2}
\begingroup
\small
\setstretch{1}
\tableofcontents
\endgroup
\clearpage

\section{Introduction and main theorem}\label{sec:intro}

Let $f:X\to Y$ be a projective morphism between smooth quasi-projective
schemes of pure dimension over a field $k$. Put
$d=\dim X-\dim Y$ and
\[
 T_f=[T_X]-f^*[T_Y]\in K_0(X).
\]
We write $f_!$ for the alternating derived direct image on $K_0$ of vector
bundles, using the identification $K_0=G_0$ on smooth schemes, and $f_*$ for
proper pushforward on Chow groups. Fix $n\geq0$ and set $m=n+d$.
We assume throughout that $m\geq0$; the case $m<0$ is trivial and is
recalled in \cref{rem:mnegative}. The degree-$n$ part of
Grothendieck--Riemann--Roch is
\begin{equation}\label{eq:rational}
 \ch_n(f_!x)=f_*\bigl[\ch(x)\Td(T_f)\bigr]_m
 \qquad\text{in }\CH^n(Y)\otimes\Q.
\end{equation}
See, for instance, \cite[Theorem~15.2]{Fulton} or \cite{SGA6}.
Here $[\,\cdot\,]_j$ denotes the homogeneous component of codimension
$j$. We use the integral Newton operations $s_j$, defined as the power sums
in the Chern roots; thus $s_0=\rk$ and $s_j=j!\ch_j$.

For $j\geq0$, let
\[
 \tau_j=\prod_{p\text{ prime}}p^{\lfloor j/(p-1)\rfloor}
\]
be the $j$-th Todd number, and put
\[
 D=\lcm(n!,\tau_m).
\]
The integer $\tau_j$ clears the denominators of the Todd polynomial in
codimension at most $j$; see \cite[Lemma~1.7.3]{Hirzebruch}. The divisibility
needed below, namely $i!\tau_{m-i}\mid\tau_m$, is proved in
\cref{sec:filtration}. In particular, the expression
\[
 P_m(x,V):=D\bigl[\ch(x)\Td(V)\bigr]_m
\]
is an integral universal polynomial in the Chern classes of $x$ and $V$.
We define the additive defect
\begin{equation}\label{eq:defect}
 \Delta_f(x)=\frac{D}{n!}s_n(f_!x)-f_*P_m(x,T_f)
 \in\CH^n(Y).
\end{equation}
By \eqref{eq:rational}, every value of $\Delta_f$ is torsion. The question
is to find a uniform annihilator as $f$, $X$, $Y$ and $x$ vary while $n$
and $m$ remain fixed.

Before stating our bound, let us recall the cases in which an integral
formula is already known.
\begin{itemize}
\item \emph{Pappas's theorem in characteristic zero.}
Pappas proves two complementary formulas, according to the sign of the
relative dimension. If $d\geq0$, his theorem gives
$\Delta_f(x)=0$ for every $x\in K_0(X)$; see
\cite[Corollary~2.3]{Pappas} and \cref{subsec:pappas}. If $d<0$
and $m=n+d\geq0$, his Theorem~2.2(b), rewritten in terms of the
relative tangent bundle and our normalization of the defect, gives
\[
 \frac{\tau_n}{D}\,\Delta_f(x)=0.
\]
We spell out this reformulation in Proposition~\ref{prop:pappas-bound}
of \cref{subsec:pappas}. Thus the two relative-dimension ranges of
Pappas's integral theorem are both accounted for in our notation.
\item \emph{Closed immersions and projective bundle projections.}
In arbitrary characteristic, the defect itself vanishes for each of
these two types of morphism:
\[
 \Delta_\iota(x)=0
 \qquad\text{and}\qquad
 \Delta_q(x)=0,
\]
where $\iota$ is a closed immersion and $q$ is a projective bundle
projection. These are the classical building blocks of the
factorization argument; see
\cite[Theorem~4.3 and Proposition~4.4]{Pappas},
\cite[Theorem~4.2]{Pappas}, and the discussion in
\cref{subsec:factorization}. Consequently, if $f$ factors as a closed
immersion into a projective bundle followed by its projection, one
obtains the classical factorization bound of
\cref{prop:factorization-bound}. Although the two constituent defects
vanish, their composition need not have zero defect: in general the
factorization estimate contains an additional Todd-number factor.
\item \emph{Finite étale morphisms.}
There is a sharper special result for finite étale covers, in the
characteristic-zero setting considered by Pappas. In this case $d=0$,
$T_f=0$, and $D=\tau_n$, so
\[
 \Delta_f(x)=\frac{\tau_n}{n!}
       \bigl(s_n(f_!x)-f_*s_n(x)\bigr).
\]
Fulton--MacPherson's covering-map formula gives
\[
 L_n\bigl(s_n(f_!x)-f_*s_n(x)\bigr)=0,
 \qquad
 L_n=\prod_{p\mid \tau_n/n!}p
     =\operatorname{rad}(\tau_n/n!);
\]
see \cite[Theorem~23.3]{FultonMacPherson}. Pappas explains in
\cite[\S2.d.1]{Pappas} that $L_n$ is the smallest universal
annihilator for this unscaled difference. The optimality argument uses
Evens--Kahn's result for classifying-space covers \cite{EvensKahn}
and Totaro's algebraic approximations \cite{Totaro}, which realize
the relevant covers by finite étale morphisms. In particular,
\[
 L_n\mid \frac{\tau_n}{n!}\mid D,
\]
and the first divisibility is often strict (for example,
$L_2=6<12=D$ and $L_3=2<24=D$). It follows that
$\Delta_f(x)=0$ exactly in this characteristic-zero setting.
This exact vanishing is also the finite étale case of Pappas's
Corollary~2.3; the additional information in the covering-map
formula is the optimal universal multiplier for the unscaled
difference. This special case is distinct from the general
projective-factorization estimate of \cref{subsec:factorization}.
\end{itemize}

We now define the integers in the main result. For $j\geq1$, set
\[
 N_j=
 \begin{cases}
  2,&j\text{ odd},\\[2pt]
  2\displaystyle\prod_{p-1\mid j}p^{\,1+v_p(j)},&j\text{ even},
 \end{cases}
 \qquad
 V_r=\prod_{j=1}^rN_j,\quad V_0=1,
\]
where the product in the definition of $N_j$ is over all primes, including
$p=2$.\footnote{For even $j\geq2$, $N_j$ is the denominator in lowest
terms of $B_j/(2j)$, with the usual Bernoulli convention $B_2=1/6$;
see \cite[Theorem~2.6]{AdamsJ}. It is also the order of the image of
the stable $J$-homomorphism
$J:\pi_{2j-1}(\mathrm{SO})\to\pi_{2j-1}^{S}$, by Adams's calculation
\cite[Theorems~1.5--1.6]{AdamsJIV} together with Quillen's proof of
the Adams conjecture \cite[Theorem~1.1]{QuillenAdams}.}
We also put $N_0=1$.

Let $F^iK_0(X)$ denote the topological filtration, generated by coherent
sheaves supported in codimension at least $i$; we set $F^i=K_0(X)$ for
$i\leq0$. The filtration is recalled in detail in \cref{sec:filtration}.

\begin{theorem}\label{thm:main}
For every $i$ with $0\leq i\leq m$ and every
$x\in F^iK_0(X)$, one has
\begin{equation}\label{eq:filtered-main}
 V_{m-i}\,\Delta_f(x)=0\qquad\text{in }\CH^n(Y).
\end{equation}
In particular, $V_m\Delta_f(x)=0$ for every $x\in K_0(X)$, or equivalently
\begin{equation}\label{eq:main}
 \frac{DV_m}{n!}s_n(f_!x)
   =f_*\bigl(V_mP_m(x,T_f)\bigr)
 \qquad\text{in }\CH^n(Y).
\end{equation}
\end{theorem}

The proof starts at the last relevant step of the topological
filtration. If $x\in F^mK_0(X)$ and $u\in\CH^m(X)$ represents its
class in $\operatorname{gr}^m_FK_0(X)$, all the terms of degree less
than $m$ vanish and
$P_m(x,T_f)=Du$. Proper pushforward respects the leading cycle, so the
target term in the defect is also $Df_*u$. Thus $\Delta_f$ vanishes
on $F^mK_0(X)$. This agreement of leading terms is the starting point
for the whole argument.

The second ingredient is a normalized Adams operation $B_a$, defined
after inverting an integer $a\geq2$. On the $i$-th graded piece of the
filtration it acts as multiplication by $a^i$, whereas algebraic
Adams--Riemann--Roch gives the defect weight $m$:
\[
 B_ax-a^ix\in F^{i+1}K_0(X)\otimes\Z[1/a],
 \qquad
 \Delta_f(B_ax)=a^m\Delta_f(x).
\]
Here the defect is extended linearly to the localized groups. These
two identities let us descend through the filtration: knowing an
annihilator on $F^{i+1}$ gives one on $F^i$ after adjoining the factor
$a^{m-i}-1$. Starting from the vanishing on $F^m$, we obtain
$\prod_{j=1}^{m-i}(a^j-1)$ as an annihilator over $\Z[1/a]$.

The last step recovers an integral identity by varying the Adams
parameter. Inverting $a$ can lose torsion information only at primes
dividing $a$. For each prime $p$, we therefore choose $a$ prime to $p$,
so that the identity over $\Z[1/a]$ also holds over $\Z_{(p)}$.
The elementary calculation in \cref{lem:Nj} lets us choose $a$ so that
the product has exactly the same $p$-adic valuation as $V_{m-i}$.
It is therefore $V_{m-i}$ times a unit in $\Z_{(p)}$, and
$V_{m-i}\Delta_f(x)=0$ in the $p$-local Chow group. Since this holds
at every prime, it holds integrally: an element of an abelian group
that vanishes in every prime localization is zero. Thus the information
lost by inverting one integer is recovered by using other Adams
operations.

The argument takes place in algebraic $K_0$ and Chow groups. Its tools
are Chern classes, the topological filtration, Bott's cannibalistic
classes and algebraic Adams--Riemann--Roch
\cite{SGA6,FultonLang,Koeck,NN}; no spectra or motivic homotopy theory
are involved. We first put the result in context in
\cref{sec:comparisons}. The proof itself occupies the next three
sections: \cref{sec:filtration} establishes the initial vanishing,
\cref{sec:adams} proves the Adams-weight identities, and
\cref{sec:annihilator} completes the induction. The article is organized to keep the proof of the main bound separate
from the classical comparisons: Section~2 provides the context, and
Sections~3--5 contain the proof.

\section{Comparison with classical integral formulas}\label{sec:comparisons}

We compare the bound of \cref{thm:main} with two sources of integral
Riemann--Roch formulas. Classical factorization gives an annihilator
in every characteristic, depending on a chosen projective embedding.
Pappas's theorem provides another comparison in characteristic zero. In
both comparisons we keep the normalization $D=\lcm(n!,\tau_m)$, so
that the integers below all annihilate the same class $\Delta_f(x)$.
Throughout this section $m\geq0$. The proof of the main theorem begins
in \cref{sec:filtration}.

\subsection{Bounds from classical factorization}\label{subsec:factorization}

Every projective morphism admits a factorization
\[
 X\xrightarrow{\,\iota\,}\PP^N_Y\xrightarrow{\,q\,}Y,
\]
where $\iota$ is a closed immersion and $q$ is the projection. The
integral Riemann--Roch formulas for these two types of morphism are
classical. Pappas proves the corresponding statements in
\cite[Theorem~4.2 and Theorem~4.3]{Pappas}; the immersion formula is
also closely related to Jouanolou's Riemann--Roch theorem without
denominators \cite{Jouanolou}. Although Pappas works throughout in
characteristic zero, the proofs of these two particular results use no
characteristic-zero hypothesis: they rely on the projective bundle
formula and on deformation to the normal cone/Koszul resolution.
We recall the statements in the form needed here, without repeating
their proofs. Their composition yields an explicit bound for the defect.

For an integer $E$ divisible by $\tau_r$, put
\[
 \Pcal_r^{(E)}(x,V)=E[\ch(x)\Td(V)]_r.
\]
This is an integral universal polynomial, by \cref{lem:todd-divisibility}
and the fact that $r!\mid\tau_r$.

\begin{lemma}\label{lem:immersion-polynomial}
Let $\iota:X\hookrightarrow P$ be a closed immersion of smooth
quasi-projective schemes of codimension $c$, with normal bundle $N_\iota$.
For every $r\geq c$, every virtual bundle $W$ on $P$, and every
multiple $E$ of $\tau_r$, one has
\begin{equation}\label{eq:immersion-polynomial}
 \Pcal_r^{(E)}(\iota_!x,W)
  =\iota_*\Pcal_{r-c}^{(E)}
       (x,\iota^*W-N_\iota).
\end{equation}
In particular, the integral Riemann--Roch defect of a closed immersion
vanishes in every degree.
\end{lemma}

\begin{lemma}\label{lem:projective-projection}
Let $q:\PP^N_Y\to Y$ be the projective-space projection. For every
$n\geq0$ and every $x\in K_0(\PP^N_Y)$, the integral formula
\[
 \frac{\tau_{n+N}}{n!}s_n(q_!x)
   =q_*\Pcal_{n+N}^{(\tau_{n+N})}(x,T_q)
\]
holds in $\CH^n(Y)$.
\end{lemma}

Pappas already points out that the classical factorization strategy
introduces additional denominators; see the Introduction and
\cite[Remark~2.6]{Pappas}. The proposition below makes this observation
quantitative: it gives an explicit annihilator arising from the chosen
factorization. This formulation is useful for comparing the classical
strategy with the uniform Adams-operation bound of \cref{thm:main}.

\begin{proposition}\label{prop:factorization-bound}
Suppose $m=n+d\geq0$ and choose a factorization
$f=q\circ\iota$ with $\iota:X\hookrightarrow\PP^N_Y$ a closed immersion.
Then
\begin{equation}\label{eq:factorization-bound}
 \frac{\tau_{n+N}}{D}\,\Delta_f(x)=0
 \qquad\text{for every }x\in K_0(X).
\end{equation}
Consequently,
\[
 \gcd\left(V_m,\frac{\tau_{n+N}}{D}\right)\Delta_f(x)=0.
\]
\end{proposition}
\begin{proof}[Sketch of proof]
Put $E=\tau_{n+N}$, so that $D\mid E$. Apply the preceding two lemmas
and use $T_f=\iota^*T_q-N_\iota$ to obtain
\[
 \frac{E}{n!}s_n(f_!x)
  =q_*\Pcal_{n+N}^{(E)}(\iota_!x,T_q)
  =f_*\Pcal_m^{(E)}(x,T_f).
\]
Since $\Pcal_m^{(E)}=(E/D)P_m$, this is
\eqref{eq:factorization-bound}. Combining it with \cref{thm:main}
gives the gcd bound.
\end{proof}

If $d\geq0$, then $D=\tau_m$. Writing $c=N-d$ for the codimension of
$\iota$, the factor in \eqref{eq:factorization-bound} becomes
\[
 \frac{\tau_{m+c}}{\tau_m}.
\]
If the chosen immersion has codimension $c=1$, this quotient is
\[
 \frac{\tau_{m+1}}{\tau_m}
   =\prod_{p-1\mid m+1}p.
\]
For example, when $m=2$ this gives the factor $2$, whereas the uniform
Adams bound is $V_2=48$. This concerns a factorization through a
codimension-one immersion; for a closed immersion itself the defect
vanishes, by \cref{lem:immersion-polynomial}. The classical bound can
therefore be much smaller than $V_m$, although it depends on the chosen
embedding. The bound in \cref{thm:main} depends only on $m$ and also
improves with the filtration degree of $x$.

\subsection{Pappas's theorem}\label{subsec:pappas}

In characteristic zero we can also use Pappas's theorem. For a smooth
scheme $Z$, a class $x\in K_0(Z)$ and a virtual bundle $V$ on $Z$, write
\[
 \CT_r(x;V)=\tau_r[\ch(x)\Td(V)]_r.
\]
This is the integral Chern--Todd polynomial. In characteristic zero,
Pappas proves that $\Delta_f(x)=0$ when $d\geq0$
\cite[Corollary~2.3]{Pappas}. In this range $D=\tau_m$, so his
formula is the exact degree-$n$ identity with our denominator-clearing
integer. His proof uses resolution of singularities and weak
factorization. The results of this paper apply over an arbitrary field.

For $d<0$ and $m=n+d\geq0$, Pappas's theorem is stated in absolute
form, with $T_X$ rather than $T_f$ on the source. More precisely, the
form we use is
\begin{equation}\label{eq:pappas-absolute}
 \CT_k(f_!x;T_Y)
   =\frac{\tau_k}{\tau_{k+d}}
       f_*\CT_{k+d}(x;T_X),
 \qquad k\geq-d,
\end{equation}
which is the content of \cite[Theorem~2.2(b)]{Pappas}. The relative
formula for the defect is obtained by combining these identities in
all degrees, not by replacing $T_X$ with $T_f$ in
\eqref{eq:pappas-absolute}.

The next proposition is not a new Riemann--Roch theorem: it is a
reformulation of Pappas's Theorem~2.2(b) in terms of the relative
tangent bundle $T_f=T_X-f^*T_Y$ and our normalization of the defect.
We include the calculation to make the resulting annihilator
$\tau_n/D$ explicit and to compare it directly with the bound
$V_m$ of \cref{thm:main}.

\begin{proposition}\label{prop:pappas-bound}
Assume $d<0$, $m=n+d\geq0$, and that the base field has characteristic
zero. Then
\[
 \frac{\tau_n}{D}\,\Delta_f(x)=0
 \qquad\text{for every }x\in K_0(X).
\]
Consequently,
\begin{equation}\label{eq:pappas-gcd}
 \gcd\left(V_m,\frac{\tau_n}{D}\right)\Delta_f(x)=0.
\end{equation}
\end{proposition}
\begin{proof}
Since $m<n$, one has $D\mid\tau_n$. Write
$t_j(V)=\tau_j[\Td(V)]_j$. For any virtual bundle $W$ and
$y\in K_0(Y)$, the identity
\[
 \ch(y)=\ch(y)\Td(W)\Td(-W)
\]
gives the universal polynomial formula
\begin{equation}\label{eq:CT-expansion}
 \frac{\tau_n}{n!}s_n(y)
  =\sum_{k=0}^n
     \frac{\tau_n}{\tau_k\tau_{n-k}}
       \CT_k(y;W)t_{n-k}(-W).
\end{equation}
All coefficients are integers by
$\tau_k\tau_{n-k}\mid\tau_n$, proved in \cref{lem:todd-divisibility}.
Apply \eqref{eq:CT-expansion} to $y=f_!x$ and $W=T_Y$. Since
$f_!x\in F^{-d}K_0(Y)$, the terms with $k<-d$ vanish. For the other
terms use \eqref{eq:pappas-absolute}, the projection formula, and put
$j=k+d$. The resulting sum is
\[
 f_*\sum_{j=0}^m
   \frac{\tau_n}{\tau_j\tau_{m-j}}
     \CT_j(x;T_X)t_{m-j}(-f^*T_Y)
 =\frac{\tau_n}{D}f_*P_m(x,T_f).
\]
The equality on the right is the same universal polynomial identity
\eqref{eq:CT-expansion}, now applied to $T_f=T_X-f^*T_Y$ and in degree
$m$. Thus
\[
 \frac{\tau_n}{n!}s_n(f_!x)
   =\frac{\tau_n}{D}f_*P_m(x,T_f),
\]
which is equivalent to $(\tau_n/D)\Delta_f(x)=0$. Combining this with
\cref{thm:main} gives \eqref{eq:pappas-gcd}.
\end{proof}

For illustration, the following table records the two annihilators in
some cases with $d<0$. The last column is the bound obtained by combining
them.
\[
\begin{array}{c|r|r|r}
(n,m)&V_m&\tau_n/D&\gcd(V_m,\tau_n/D)\\ \hline
(3,2)&48&2&2\\
(4,2)&48&30&6\\
(5,2)&48&12&12\\
(6,2)&48&84&12\\
(8,3)&96&90&6
\end{array}
\]
Thus neither annihilator is uniformly smaller: Pappas gives $2$
instead of $48$ at $(n,m)=(3,2)$, whereas $V_m$ gives $48$ instead of
$84$ at $(6,2)$. In the latter case their gcd gives the still smaller
factor $12$. The two bounds can therefore be used together. Thus the
content of Proposition~\ref{prop:pappas-bound} is the explicit
translation of Pappas's formula into our relative-defect notation;
the additional comparison is the combination with the independent
Adams bound of Theorem~\ref{thm:main}.

\section{The topological filtration and the initial vanishing}\label{sec:filtration}

We begin with two standard properties of the topological filtration.
The first identifies the leading Chern class of a class supported in a
fixed codimension. The second describes proper pushforward on the
associated graded group. Together they show that the defect is already
zero on the last relevant filtration step.

For a virtual bundle $V$, write
\[
 t_j(V)=\tau_j[\Td(V)]_j.
\]
This is an integral polynomial in the Chern classes of $V$.

\begin{lemma}\label{lem:todd-divisibility}
For $0\leq i\leq m$ one has
\[
 i!\tau_{m-i}\mid\tau_m
 \qquad\text{and}\qquad
 \tau_i\tau_{m-i}\mid\tau_m.
\]
\end{lemma}
These divisibility properties are particular cases of \cite[Lemma~2.1]{Pappas}. The case $m=0$ is immediate.

It follows that the integral polynomial in \eqref{eq:defect} can be
written as
\begin{equation}\label{eq:Pm}
 P_m(x,V)=\sum_{i=0}^m
       \frac{D}{i!\tau_{m-i}}\,s_i(x)t_{m-i}(V).
\end{equation}
In particular, it is additive in $x$ and is defined in integral Chow
groups, including when those groups have torsion.

Let $F^iK_0(X)$ be the subgroup generated by classes of coherent sheaves
whose supports have codimension at least $i$. Since $X$ is regular,
$G_0(X)\simeq K_0(X)$, so this defines a filtration of $K_0(X)$. We put
$F^i=K_0(X)$ for $i\leq0$ and $F^i=0$ for $i>\dim X$. Its associated
graded group is denoted by $\operatorname{gr}^i_FK_0(X)=F^i/F^{i+1}$.
The cycle map
\[
 \cl_i:\CH^i(X)\longrightarrow\operatorname{gr}^i_FK_0(X),
 \qquad [Z]\longmapsto[\OO_Z],
\]
is surjective.

\begin{lemma}\label{lem:leading}
Let $x\in F^iK_0(X)$ and let $u\in\CH^i(X)$ represent its class in
$\operatorname{gr}^i_FK_0(X)$. Then
\begin{equation}\label{eq:leading}
 s_j(x)=0\quad(0\leq j<i),\qquad s_i(x)=i!u.
\end{equation}
\end{lemma}

This is a standard consequence of Riemann--Roch without denominators
and Newton's identities; see \cite[Example~15.3.6]{Fulton} and
\cite{Jouanolou}.

\begin{lemma}\label{lem:pushforward}
The proper pushforward satisfies
\[
 f_!(F^iK_0(X))\subset F^{i-d}K_0(Y).
\]
If $i-d\geq0$, the following square commutes:
\begin{equation}\label{eq:cycle}
\begin{tikzcd}[column sep=4.5em,row sep=3em]
 \CH^i(X) \arrow[r,"\cl_i"] \arrow[d,"f_*"'] &
 \operatorname{gr}_F^iK_0(X) \arrow[d,"\operatorname{gr}(f_!)"] \\
 \CH^{i-d}(Y) \arrow[r,"\cl_{i-d}"'] &
 \operatorname{gr}_F^{i-d}K_0(Y).
\end{tikzcd}
\end{equation}
\end{lemma}

For the compatibility of the cycle map with proper pushforward,
see \cite[\S\S3.1--3.2]{HM}.

\Needspace{12\baselineskip}
\begin{lemma}\label{lem:top}
For $m\geq0$, the defect $\Delta_f$ vanishes on $F^mK_0(X)$.
\end{lemma}
\begin{proof}
Take $x\in F^mK_0(X)$ and choose $u\in\CH^m(X)$ representing its
class in $\operatorname{gr}^m_FK_0(X)$. By \cref{lem:leading} and
\eqref{eq:Pm}, only the term $i=m$ contributes to $P_m(x,T_f)$; it is
$Du$. By \cref{lem:pushforward}, the class of $f_!x$ in
$\operatorname{gr}^n_FK_0(Y)$ is represented by $f_*u$. Thus
\[
 \frac{D}{n!}s_n(f_!x)=Df_*u=f_*P_m(x,T_f),
\]
which is the required vanishing.
\end{proof}

\begin{remark}\label{rem:mnegative}
If $m<0$, then $-d>n$ and the support estimate gives
$f_!K_0(X)\subset F^{-d}K_0(Y)\subset F^{n+1}K_0(Y)$. Hence
$s_n(f_!x)=0$ integrally for every $x\in K_0(X)$.
\end{remark}

\section{Adams operations and their weights}\label{sec:adams}

We next construct an operation which is scalar on each associated-graded
piece of the topological filtration. Its definition involves the Bott
cannibalistic class of the virtual tangent bundle. We then check separately
its effect on the target and source terms of the defect. Keeping these
three facts separate makes the weight calculation transparent.

Fix an integer $a\geq2$ and put $R=\Z[1/a]$. We write
$K_0(X)_R=K_0(X)\otimes R$ and
$\CH^j(X)_R=\CH^j(X)\otimes R$. The defect extends $R$-linearly to
\[
 \Delta_{f,R}:K_0(X)_R\longrightarrow\CH^n(Y)_R.
\]
The Adams operation $\psi^a$ is characterized by
$\psi^a[L]=[L^{\otimes a}]$ for line bundles. The Bott class $\theta^a$
is multiplicative for direct sums and satisfies
\[
 \theta^a(L)=1+[L]+\cdots+[L]^{a-1}
\]
for a line bundle $L$. For a virtual bundle it is defined by
multiplicativity.

If $E$ has rank $e$, then $\theta^a(E)$ has rank $a^e$. The ideal of
rank-zero classes in $K_0(X)$ is nilpotent: it is contained in $F^1$,
and the topological filtration is multiplicative and finite. Consequently,
a class whose rank is a unit in $R$ is a unit in $K_0(X)_R$. In
particular, $\theta^a(E)$ is invertible in $K_0(X)_R$ for every virtual
bundle $E$.

\Needspace{8\baselineskip}
Write $T_f^\vee=[\Omega_X]-f^*[\Omega_Y]$ for the dual virtual
bundle. Define the normalized Adams operation
\begin{equation}\label{eq:B}
 B_a(x)=a^d\,\theta^a(T_f^\vee)^{-1}\,\psi^a(x)
 \qquad (x\in K_0(X)_R).
\end{equation}
The factor multiplying $\psi^a$ has rank one.

\begin{lemma}\label{lem:adams-weight}
For $x\in F^iK_0(X)_R$ one has
\[
 B_a(x)-a^i x\in F^{i+1}K_0(X)_R.
\]
\end{lemma}
\begin{proof}
By \cite[Definition~4.13 and Proposition~5.1]{HM}, the homological Adams
operation $\theta^a(T_X^\vee)^{-1}\psi^a$ preserves the topological
filtration over $R$ and has weight $a^{i-\dim X}$ in codimension $i$.
Multiplication by a class acts on each graded piece by its rank. Since
$\theta^a(T_X^\vee)$ has rank $a^{\dim X}$, the operation $\psi^a$ has
weight $a^i$. The rank-one factor in \eqref{eq:B} does not change this
weight, as required.
\end{proof}

\begin{lemma}\label{lem:adams-target}
For every $x\in K_0(X)_R$,
\[
 \frac{D}{n!}s_n(f_!B_ax)=a^m\frac{D}{n!}s_n(f_!x).
\]
\end{lemma}
\begin{proof}
The algebraic Adams--Riemann--Roch formula, in the form needed here, is
\begin{equation}\label{eq:ARR}
 \psi^a(f_!x)=f_!\bigl(\theta^a(T_f^\vee)^{-1}\psi^a(x)\bigr)
 \quad\text{in }K_0(Y)_R.
\end{equation}
It is the usual formula for projective local complete intersection
morphisms; in our situation the virtual tangent bundle is $T_f$.
See \cite[Exp.~V]{SGA6}, \cite[Ch.~V]{FultonLang},
\cite{Koeck} and \cite[Corollary~4.5]{NN}.
By \eqref{eq:B} and \eqref{eq:ARR},
\[
 f_!(B_ax)=a^d\psi^a(f_!x).
\]
The identity $s_n(\psi^a y)=a^n s_n(y)$ is an integral universal
identity: it follows immediately after writing the Chern roots of $y$
as line roots. The displayed equality now follows from $m=n+d$.
\end{proof}

\begin{lemma}\label{lem:universal-detection}
Let $R$ be a localization of $\Z$. An identity between universal
polynomials with coefficients in $R$ in the Chern classes of finitely
many vector bundles, if valid with rational coefficients, is valid in
$\CH^*(X)\otimes R$ for every smooth quasi-projective scheme $X$ over a
field.
\end{lemma}
\begin{proof}
Pull back to the product of the flag bundles of the finitely many
bundles. Pullback on Chow groups is injective, and all the bundles split
there into sums of line bundles. On this flag bundle, each line bundle
can be written as $A\otimes C^{-1}$ with $A$ and $C$ globally generated:
twist by a sufficiently high power of a very ample line bundle. The
sections of the finitely many $A$'s and $C$'s define a morphism to a
finite product of projective spaces, and the line bundles are pullbacks
of tensor products of the hyperplane line bundles and their inverses.
The Chow ring of this product of projective spaces is free over $\Z$.
Thus a polynomial identity which holds rationally holds there over
$R$. Pulling back and using injectivity proves the assertion.
\end{proof}

\begin{lemma}\label{lem:adams-source}
For every $x\in K_0(X)_R$,
\[
 P_m(B_ax,T_f)=a^mP_m(x,T_f)
 \quad\text{in }\CH^m(X)_R.
\]
\end{lemma}
\begin{proof}
Let $\Phi^a$ act on $\CH^j(X)\otimes\Q$ as multiplication by $a^j$.
For a line bundle $L$ with first Chern class $t$, the identity
\[
 \frac{t}{1-e^{-t}}\,
 \frac{a}{1+e^{-t}+\cdots+e^{-(a-1)t}}
 =\frac{at}{1-e^{-at}}
\]
shows that
\[
 \Td(V)\,\ch\bigl(a^e\theta^a(V^\vee)^{-1}\bigr)
       =\Phi^a\Td(V)
\]
for every virtual bundle $V$ of rank $e$, by multiplicativity. Since
$\ch(\psi^a x)=\Phi^a\ch(x)$, taking the component of degree $m$ gives
\[
 P_m\bigl(a^e\theta^a(V^\vee)^{-1}\psi^a(x),V\bigr)
       =a^mP_m(x,V)
\]
with rational coefficients. Both sides are universal polynomials over
$R$, so \cref{lem:universal-detection} makes the identity integral over
$R$. Take $V=T_f$, whose rank is $d$, to obtain the assertion.
\end{proof}

\begin{proposition}\label{prop:weight-defect}
For every $x\in K_0(X)_R$ one has
\begin{equation}\label{eq:defect-weight}
 \Delta_{f,R}(B_ax)=a^m\Delta_{f,R}(x).
\end{equation}
\end{proposition}
\begin{proof}
The first term in the definition of the defect is multiplied by $a^m$
by \cref{lem:adams-target}; the second is multiplied by $a^m$ by
\cref{lem:adams-source} and the projection formula. Subtracting gives
\eqref{eq:defect-weight}.
\end{proof}

\section{The annihilator}\label{sec:annihilator}

We now apply the weight identity to the filtration. The induction is
short: if the estimate is known on $F^{i+1}$, the difference between
$B_a x$ and $a^i x$ lies in that subgroup, while the defect of $B_a x$
is known exactly. This adds the factor $a^{m-i}-1$ at the $i$-th step.

\begin{proposition}\label{prop:local}
Let $0\leq i\leq m$ and $x\in F^iK_0(X)$. For every integer $a\geq2$,
\begin{equation}\label{eq:local}
 \left(\prod_{j=1}^{m-i}(a^j-1)\right)\Delta_f(x)=0
 \qquad\text{in }\CH^n(Y)\otimes\Z[1/a].
\end{equation}
\end{proposition}
\begin{proof}
We argue by descending induction on $i$. For $i=m$, the product is
empty and the assertion is \cref{lem:top}. Suppose the result is known
for $i+1$, and take $x\in F^iK_0(X)$. By \cref{lem:adams-weight},
\[
 x'=B_ax-a^ix\in F^{i+1}K_0(X)_R.
\]
The induction hypothesis applies to $x'$ by $R$-linearity. On the other
hand, \cref{prop:weight-defect} gives
\[
 \Delta_{f,R}(x')
   =\Delta_{f,R}(B_ax)-a^i\Delta_{f,R}(x)
   =a^i(a^{m-i}-1)\Delta_{f,R}(x).
\]
Since $a^i$ is a unit in $R$, the new factor is $a^{m-i}-1$. This
completes the induction.
\end{proof}

The local estimate becomes an integral one by choosing the Adams
parameter separately at each prime. The integers that occur in this
choice are classical.

\Needspace{12\baselineskip}
\begin{lemma}\label{lem:Nj}
For every prime $p$ and every $j\geq1$,
\[
 v_p(N_j)=\min_{\substack{a\geq2\\p\nmid a}}v_p(a^j-1).
\]
In particular, $N_j=2$ for odd $j$, while for even $j$ it is the
denominator of the Bernoulli number $B_j/(2j)$. These are the classical Adams integers; for even $j$ they are the
denominators that occur in the image-of-$J$ calculation of Adams
\cite{AdamsJ}; see also \cite[\S5.1]{HM}.
\end{lemma}
\begin{proof}
For an odd prime $p$, if $p-1\nmid j$, there is a unit $a$ modulo $p$
with $a^j\not\equiv1\pmod p$, so the minimum valuation is zero. If
$p-1\mid j$, every unit satisfies $a^j\equiv1\pmod p$, and choosing a
generator of $(\Z/p^2)^\times$ gives the minimum $1+v_p(j)$ by the
lifting-the-exponent formula. For $p=2$, take odd $a$. If $j$ is odd,
the minimum of $v_2(a^j-1)$ is $1$, attained at $a=3$. If $j$ is even,
the minimum is $2+v_2(j)$, again attained at $a=3$. These are precisely
the valuations in the definition of $N_j$. The Bernoulli denominator
identity for even $j$ is the classical von Staudt--Clausen formula.
\end{proof}

\begin{proof}[Proof of \cref{thm:main}]
Fix $i$ with $0\leq i\leq m$ and $x\in F^iK_0(X)$. For each prime
$p$, if $p$ is odd choose $a\geq2$ whose residue class generates
$(\Z/p^2)^\times$; if $p=2$, choose $a=3$. The valuation calculation
in \cref{lem:Nj} gives
\[
 v_p\left(\prod_{j=1}^{m-i}(a^j-1)\right)
   =\sum_{j=1}^{m-i}v_p(N_j)=v_p(V_{m-i}).
\]
Thus the product in \cref{eq:local} differs from $V_{m-i}$ by a unit
in $\Z_{(p)}$. It follows that $V_{m-i}\Delta_f(x)=0$ after localizing
at every prime. An element of an abelian group which vanishes in every
prime localization is zero. This proves \eqref{eq:filtered-main}
for every $i$, and hence the theorem.
\end{proof}

\begin{corollary}\label{cor:localized}
Let $N=\max(n,m+1)$. If $m\geq0$, the denominator-cleared formula
\[
 \frac{D}{n!}s_n(f_!x)=f_*P_m(x,T_f)
\]
holds in $\CH^n(Y)\otimes\Z[1/N!]$. 
\end{corollary}
\begin{proof}
Every prime divisor of $N_j$ is at most $j+1$, so every prime divisor
of $V_m$ is at most $m+1$. Every prime divisor of $D$ is at most
$\max(n,m+1)$: this is immediate for $n!$, and for $\tau_m$ it follows
from $p-1\leq m$. Thus both $D$ and $V_m$ become units after tensoring
with $\Z[1/N!]$. The result follows from \cref{thm:main} and
\cref{rem:mnegative}.
\end{proof}

Since $D$ is a unit in $\Z[1/N!]$ when $m\geq0$, we may divide the
identity in \cref{cor:localized} by $D$ to obtain
\[
 \ch_n(f_!x)=f_*[\ch(x)\Td(T_f)]_m
 \qquad\text{in }\CH^n(Y)\otimes\Z[1/N!] \ ,
\]
that is the Riemann--Roch formula
proved in our companion paper \cite[Theorem~1.1]{NNAway}.

\section*{Use of AI}
The authors used ChatGPT (OpenAI) as an AI-assisted tool during the
development and formulation of parts of the arguments and in the
preparation and revision of the exposition. The authors take
responsibility for the mathematical content and for the final version
of the manuscript.

\section*{Funding}
Both authors were supported by FEDER/AEI/MICINN under Grant
No.~PID2022-142024NB-I00 (SIGA).

\begingroup
\footnotesize
\setstretch{1}

\endgroup

\par\medskip
\begin{flushleft}
\footnotesize
\setstretch{1}
\paperemails
\end{flushleft}


\begin{thebibliography}{99}
\setlength{\itemsep}{0pt}
\setlength{\parsep}{0pt}
\setlength{\parskip}{0pt}

\bibitem{AdamsJ}
J.~F.~Adams,
\href{https://doi.org/10.1016/0040-9383(65)90040-6}{On the groups $J(X)$--II},
\textit{Topology} \textbf{3} (1965), 137--171.


\bibitem{AdamsJIV}
J.~F.~Adams,
\href{https://www.sas.rochester.edu/mth/sites/doug-ravenel/otherpapers/Adams-JX4.pdf}{On the groups $J(X)$--IV},
\textit{Topology} \textbf{5} (1966), 21--71.


\bibitem{EvensKahn}
L.~Evens and D.~Kahn,
An integral Riemann--Roch formula for induced representations of finite groups,
\textit{Trans. Amer. Math. Soc.} \textbf{245} (1978), 331--347.


\bibitem{Fulton}
W.~Fulton,
\textit{Intersection Theory}, 2nd ed.,
Ergebnisse der Mathematik und ihrer Grenzgebiete (3), vol.~2,
Springer-Verlag, Berlin, 1998.


\bibitem{FultonMacPherson}
W.~Fulton and R.~MacPherson,
\href{https://doi.org/10.2307/1971288}{Characteristic classes of direct image bundles for covering maps},
\textit{Ann. of Math. (2)} \textbf{125} (1987), no.~1, 1--92.


\bibitem{FultonLang}
W.~Fulton and S.~Lang,
\textit{Riemann--Roch Algebra},
Grundlehren der mathematischen Wissenschaften, vol.~277,
Springer-Verlag, New York, 1985.


\bibitem{Hirzebruch}
F.~Hirzebruch,
\textit{Topological Methods in Algebraic Geometry},
reprint of the 1978 edition, Classics in Mathematics,
Springer-Verlag, Berlin, 1995.


\bibitem{HM}
O.~Haution and A.~S.~Merkurjev,
\href{https://doi.org/10.4171/EMSS/50}{Connective $K$-theory and Adams operations},
\textit{EMS Surv. Math. Sci.} \textbf{8} (2021), no.~1/2, 135--162.


\bibitem{Jouanolou}
J.-P.~Jouanolou,
\textit{Riemann--Roch sans dénominateurs},
\textit{Invent. Math.} \textbf{11} (1970), 15--26.


\bibitem{Koeck}
B.~Köck,
\href{https://doi.org/10.1016/S0012-9593(98)80140-7}{The Grothendieck--Riemann--Roch theorem for group scheme actions},
\textit{Ann. Sci. École Norm. Sup.} (4) \textbf{31} (1998), no.~3,
415--458.


\bibitem{NN}
A.~Navarro and J.~Navarro,
\href{https://doi.org/10.1016/j.exmath.2023.07.002}{The Riemann--Roch theorem for the Adams operations},
\textit{Expo. Math.} \textbf{41} (2023), no.~4, article 125513.


\bibitem{NNAway}
A.~Navarro Garmendia and J.~Navarro Garmendia,
\textit{The Riemann--Roch formula away from the small primes}, manuscript, 2026.


\bibitem{Pappas}
G.~Pappas,
\href{https://doi.org/10.1007/s00222-007-0067-9}{Integral Grothendieck--Riemann--Roch theorem},
\textit{Invent. Math.} \textbf{170} (2007), no.~3, 455--481.


\bibitem{QuillenAdams}
D.~Quillen,
\href{https://math.univ-cotedazur.fr/u/cazanave/Gdt/ImJ/Quillen.pdf}{The Adams conjecture},
\textit{Topology} \textbf{10} (1971), 67--80.


\bibitem{SGA6}
P.~Berthelot, A.~Grothendieck and L.~Illusie (eds.),
\textit{Théorie des intersections et théorème de Riemann--Roch},
Séminaire de Géométrie Algébrique du Bois Marie 1966--1967 (SGA~6),
Lecture Notes in Mathematics, vol.~225, Springer-Verlag, Berlin, 1971.


\bibitem{Totaro}
B.~Totaro,
The Chow ring of a classifying space,
in \textit{Algebraic $K$-theory} (Seattle, WA, 1997),
Proc. Sympos. Pure Math., vol.~67, pp.~249--281,
American Mathematical Society, Providence, RI, 1999.


\end{thebibliography}
\end{document}